\documentclass[a4paper]{article}
\usepackage{latexsym}
\usepackage{listings}
\usepackage{amsthm,amsmath,amssymb}
\usepackage[dvips]{graphicx, color}
\usepackage[title]{appendix}
\usepackage{float}
\usepackage[linesnumbered,ruled,vlined]{algorithm2e}
\usepackage{color}
\usepackage{bbm}

\usepackage{tikz}
\usetikzlibrary{quantikz}
\usepackage{mathtools}

\newtheorem{thm}{Theorem}

\newtheorem*{thm*}{Theorem}

\newtheorem{obs}[]{Observation}

\newtheorem*{ex*}{Example}
\newtheorem{ex}[]{Example}

\newtheorem{lem}[thm]{Lemma}

\newtheorem{rem}[]{Remark}

\newcommand{\cE}{{\mathcal E}}

\newcommand{\E}{{\Bbb E}}
\newcommand{\F}{{\Bbb F}}

\newcommand{\N}{{\Bbb N}}

\newcommand{\Q}{{\Bbb Q}\hspace{.06em}}
\newcommand{\R}{{\Bbb R}}

\newcommand{\Z}{{\Bbb Z}}

\def\={\:=\:}  \def\+{\,+\,}

\def\a{\alpha} \def\b{{\beta}}  \def\ba{\overline\a}

\def\be{\begin{equation}}   \def\ee{\end{equation}}
\def\bes{\begin{equation*}}   \def\ees{\end{equation*}}
\def\ba{\begin{aligned}}   \def\ea{\end{aligned}}
\def\bc{\begin{cases}}   \def\ec{\end{cases}}
\def\bp{\begin{proof}}   \def\ep{\end{proof}}

\def\SL{\mathrm{SL}}

\def\qqan{\qquad\mathrm{and}\qquad}
\def\qan{\quad\mathrm{and}\quad}

\newcommand{\dis}{\displaystyle}

\def\SL{\mathrm{SL}}

\def\bbm1{\mathbbm 1}

\def\wh{\widehat}
\def\be{\begin{equation}}   \def\ee{\end{equation}}
\def\bes{\begin{equation*}}   \def\ees{\end{equation*}}
\def\bea{\begin{equation}\begin{aligned}}   
\def\eea{\end{aligned}\end{equation}}

\def\Pic{\mathrm{Pic}}

\def\bm{\begin{matrix}}
\def\em{\end{matrix}}
\def\bpm{\begin{pmatrix}}
\def\epm{\end{pmatrix}}

\begin{document}
\title{\bf Murmurations\\ and Sato-Tate Conjectures\\ for High Rank Zetas 
of Elliptic Curves} 
\author{{\bf Zhan SHI\footnote{Graduate Program of Mathematics for Innovation, 
Kyushu University, Fukuoka, Japan}} and 
{\bf Lin WENG\footnote{Graduate School of Mathematics, Kyushu University, Fukuoka, Japan}}}  
\date{(October 7, 2024 in Fukuoka)}
\maketitle
\begin{abstract}

For elliptic curves over rationals, there are a well-known conjecture of Sato-Tate and a new computational 
guided murmuration phenomenon, for which the abelian Hasse-Weil zeta functions are used.  In this paper, 
we show that both the murmurations and the Sato-Tate conjecture stand equally well for non-abelian high 
rank zeta functions of the $p$-reductions of elliptic curves over rationals. We establish our results by 
carefully examining asymptotic  behaviors of the $p$-reduction invariants 
$a_{E/\mathbb F_p,n}\ (n\geq 1)$, the rank $n$ analogous of the rank one $a$-invariant 
$a_{E/\mathbb F_p}=1+p-N_{E/\mathbb F_p}$ of elliptic curve $E/\mathbb F_{p}$. Such asymptotic results 
are based on the counting miracle of the so-called $\alpha_{E/\mathbb F_q,n}$- and 
$\beta_{E/\mathbb F_q,n}$-invariants of $E/\mathbb F_q$ in rank $n$, and a remarkable recursive relation 
on the $\beta_{E/\mathbb F_q,n}$-invariants, both established by Weng-Zagier in \cite{EC}.
\end{abstract}


\section{Statement of Main Theorem}
Murmuration is an intriguing phenomenon newly discovered by He-Lee-Oliver-Pozdnyakov in \cite{MUR1}
for families of elliptic curves $\E$ defined over the field $\Q$ of rationals, offering an aesthetic  intuitive 
relation between the averages of the coefficients $a_{E/\F_{p_i}}$ for the $p_i$-reductions $E/\F_{p_i}$ of 
$\E/\Q$ and the arithmetic ranks of  $\E(\Q)$. On the other hand, the Sato-Tate conjecture is a statistical 
statement formulated around 1960, exposing a natural density function for the distributions of the $a$-
invariant $a_{E/\F_{p_i}}$'s of a fixed elliptic curve $\E/\Q$, as  a secondary structure beyond the Riemann 
hypothesis.  Both of these two are  obviously closely associated to the classical Hasse-Weil zeta functions 
for elliptic curves over finite fields. While the mystery surrounding  murmurations of elliptic curves remains 
to be understood, the Sato-Tate conjecture (and its generalization to all totally real fields) was proved by 
Clozel-Harris-Shepherd(-)Barron-Taylor under mild assumptions in 2008, and completed by Barnet(-)Lamb-
Geraghty-Harris-Taylor in 2011. Several generalizations to other algebraic varieties and fields have been 
made and are widely open.

In this paper, we show that both murmurations and Sato-Tate conjecture stand well when higher rank 
non-abelian zeta functions of elliptic curves are used, instead of the classical Hasse-Weil zeta functions. 
Recall that, for an integral regular projective curve $X$ over a finite field $\F_q$, its high rank zeta functions $\wh\zeta_{X/\F_q;n}(s)$  were introduced by the senior author in \cite{HRZ2} to study rank $n$ semi-stable vector bundles over $X/\F_q$, as a genuine generalization of the classical congruent Artin-Weil zeta function $\wh\zeta_{X/\F_q}(s)$:
$$
\wh\zeta_{X/\F_q}(s)=\wh\zeta_{X/\F_q;1}(s).
$$ 
We establish our results on the high rank Sato-Tate and high rank murmurations for elliptic curves $\E/\Q$
by carefully examining asymptotic  behaviors of the invariants $a_{E/\F_{p_i},n}$ of the $p_i$-reduction $E/\F_{p_i}$ of $\E/\Q$, the rank $n$ analogous of the rank one  $a$-invariant $a_{E/\F_{p_i}}:=1+p_i-\#E(\F_{p_i})$. Here $p_i$ denotes the $i$-th prime number. As to be seen below, these asymptotic results are established based on some structural recursion relations, which themselves are obtained from the \lq counting miracle', on the so-called $\a_{E/\F_p,n}$- and $\b_{E/\F_p,n}$-invariants for elliptic curves $E/\F_p$. In fact,  the \lq counting miracle', relating $\a_{E/\F_p,n}$- and $\b_{E/\F_p,n}$-invariants,  was first proved by Weng-Zagier in \cite{EC} for elliptic curves $E/\F_p$ using Atiyah bundles and some combinatorial techniques, then generalized by Sugahara \cite{SGHR} and  Mozgovoy-Reineke independently \cite{MR} for curves over finite fields, and the structural recursion relations on the  $\b_{E/\F_p,n}$-invariants are established  in \cite{EC} through some complicated combinatorial discussions.

To facilitate our ensuing discussion, next, we give a quick review of rank $n$ zeta function for curves over finite fields following \cite{HRZ2}, see also \cite{SLn}. For a  fixed $n\in\Z_{\geq 1}$, let $X$ be a regular projective (integral) curve of genus $g$ defined over a finite filed $\F_q$ with $q$ elements. It is well known that  the classical Artin-Weil zeta function of $X/\F_q$ can be defined using the following power series
\be
\zeta_{X/\F_q}(s):=\exp\left(\sum_{k=1}^\infty\frac{\#X(\F_{q^k})}{k}(q^{-s})^k\right)\qquad\Re(s)>1.
\ee
Here, as usual, $X(\F_{q^k})$ denotes the set of $\F_{q^k}$-rational points of $X$.
In addition, with a change of  variables $t=q^{-s}$, we arrive at the 
corresponding Zeta function for $X/\F_q$: 
\be
Z_{X/\F_q}(t):=\zeta_{X/\F_q}(s)=\exp\left(\sum_{k=1}^\infty\frac{\#X(\F_{q^k})}{k}t^k\right).
\ee
The following theorem is well known:
\begin{thm}[Zeta Properties, see e.g. \cite{H}]\label{thm1}
Let $X$ be an integral regular projective curve of genus $g$ over a finite filed $\F_q$. Then
\begin{enumerate}
\item [(1)] (Rationality) 
$Z_{X/\F_q}(t)$ is a rational function of the form
$$
Z_{X/\F_q}(t)=\frac{P_{X/\F_q}(t)}{(1-t)(1-qt)}.
$$
where $P_{X/\F_q}(t)$ is a polynomial of degree $2g$ with integer coefficients. 
\item [(2)] (Functional Equation) 
$\zeta_{X/\F_q}(s)$ satisfies the standard functional equation
$$
\zeta_{X/\F_q}(1-s)=q^{(g-1)(2s-1)}\cdot \zeta_{X/\F_q}(s).
$$
\item [(3)] (Riemann Hypothesis)
If $ \zeta_{X/\F_q}(s)=0$, then $\Re(s)= \frac{1}{2}$.
\end{enumerate}
\end{thm}

It turns out that the classical Artin-Zeta function $\zeta_{X/\F_q}(s)$ above is an abelian starting point of a family of the so-called  rank $n$ zeta functions for  $X/\F_q$ (n=1,2,\ldots). Indeed, by using  the Euler product and the ramification theory on the  base changes of the constant fields $\F_q\mapsto\F_{q^k}$, 
it is not difficulty to arrive at
$$\zeta_{X/\F_q}(s)=\sum_{d\in \Z}\sum_{L\in \Pic^d(X/\F_q)}\frac{q^{h^0(X,L)}-1}{q-1}(q^{-s})^{\deg(L)}\qquad\Re(s)> 1.$$
Here $h^0(X,L):=\#H^0(X,L)$ denotes the 0th cohomology of the line bundle $L$ on $X$, and $\deg(L)$ denotes the degree of $L$.

The first version of a rank $n$ non-abelian zeta function for the curve $X/\F_q$ was  introduced by the senior author in \cite{HRZ1}, by counting \underline{all} rank $n$ semi-stable vector bundles on $X/\F_q$. Even these newly defined functions satisfy most of the zeta properties, including the functional equation and the rationality, but fatally they do not satisfy the Riemann hypothesis as examples indicate. This means that our counting at the moment was not correct. Motivated by Drinfeld's work  (\cite{D}) on counting two-dimensional irreducible representations of the fundamental group of curves over finite fields,  an essential genuine restriction on rank $n$ semi-stable vector bundles is introduced, namely, counting only those rank $n$ semi-stable vector bundles on $X/\F_q$ whose degrees are \underline{multiples of  $n$}. This then leads to the right current formulation of the {\it rank $n$ zeta function $\zeta_{X/\F_q,n}(s)$ of $X/\F_q$} in  \cite{HRZ2}:    in the region $\Re(s)>1$, 
\be
\zeta_{X/\F_q,n}(s) := \zeta_{X,n}(s) := \sum_{d\in n\Z}\sum_{V}\frac{q^{h^0(X,V)}-1}{\#\mbox{Aut}(V)}(q^{-s})^{\deg(V)},
\ee
where the second sum on $V$ rums over all  rank $n$ ($\F_q$-rational) semi-stable vector bundle on $X/\F_q$ of degree $d=mn\in n\Z$,  $\mbox{Aut}(V)$ denotes the automorphism group of the vector bundle $V$, and $h^0(X,V)$, resp. $\deg(V)$, denotes the  0th cohomology, resp. the degree of $V$.

As in \cite{HRZ2}, we introduce the following $\a$- and $\b$- invariants of $X/\F_p$:
$$
\alpha_{X/\F_q,n}(d):=\sum_{V}\frac{q^{h^0(X,V)}-1}{\#\mbox{Aut}(V)},
\qquad 
\beta_{X/\F_q,n}(d)=\sum_{V}\frac{1}{\#\mbox{Aut}(V)},
$$
where similarly, in both of the  summations,  $V$ rums over all  rank $n$ ($\F_q$-rational) semi-stable vector bundle on $X/\F_q$ of degree $d.$
Then, easily, the rank $n$ zeta function $\zeta_{X,n}(s)$ can be rewritten as a generating function of $\alpha_{X,n}(mn)$:
\be\label{eq8}
\zeta_{X/\F_q,n}(s)=\sum_{d\in n\Z_{\geq 0}}\alpha_{X/\F_q,n}(d)t^d=\sum_{m=0}^{\infty}\alpha_{X/\F_q,n}(mn)T^m,
\ee
where we have set  $T_n:=t^n=Q_n^{-s}$ with $Q_n:=q^n$. When $n$ is clear from the text, we often omit the index $n$. For examples, write $Q$ for $Q_n$ and $T$ for $T_n$.

Similar to Theorem\,\ref{thm1}, by using 
the Riemann-Roch theorem, the duality  and a vanishing theorem for semi-stable 
bundles, we have the following:
\begin{thm}\label{thm2}[Zeta Properties \cite{HRZ2}, see also \cite{SLn}] Fix a natural number $n$. The rank $n$ non-abelian zeta function 
$\zeta_{X,n}(s)$ of an integral regular projective curve $X/\F_q$ satisfies the following standard zeta properties:
\begin{enumerate}
\item (Naturality) The rank one zeta function $\zeta_{X,1}(s)$ coincides with the classical Artin-Weil zeta function 
$\zeta_{X}(s)$ of $X/\F_q$.
\item (Rationality) There exists a polynomial $P_{X,n}(T)\in \Q[T]$ of degree $2g$, such that
$$
\zeta_{X,n}(s)=\frac{P_{X,n}(T)}{(1-T)(1-QT)}.
$$
\item (Functional Equation) $\zeta_{X,n}(s)$ satisfies the standard functional equation
$$
\zeta_{X,n}(1-s)=Q^{(g-1)(2s-1)}\cdot \zeta_{X,n}(s).
$$
\end{enumerate}
\end{thm}
 Furthermore, it is conjectured in \cite{HRZ2} that    all these rank $n$ zeta functions $\zeta_{X,\F_q;n}(s)$ satisfy the Riemann Hypothesis. Surprisingly, this conjecture remains widely open, even its number theoretic analogue has been established (except when $n=1$ for  the lack of symmetry), up to a finite box depnding on $n$ (\cite{W3}). The first major breakthrough in this direction is the following:
  \begin{thm}[$n=1$: Hasse \cite{Sil}, $n\geq 2$: Weng-Zagier \cite{EC}]\label{thm3} Let $E/\F_q$ be an elliptic curve. Then
 $\zeta_{E/\F_q,n}(s)$ satisfies the Riemann hypothesis. That is to say,\\
 $\zeta_{E/\F_q,n}(s)=0$ implies that $\Re(s)=\frac{1}{2}$.\footnote{We mention in passing that, besides this elliptic curve case, the Riemann hypothesis for $\zeta_{X/\F_q,n}(s)$ has been established  successfully  when
 \begin{enumerate} 
 \item[(i)] $n=2$ by H. Yoshida, see e.g. \S2 of arXiv:2201.03703.
 \item [(ii)] $n=3$ by Weng in \lq Riemann Hypothesis for Non-Abelian Zeta Functions of Curves over Finite Fields', arXiv:2201.03703.
 \item [(iii)] $g=2$ asymptotically by Shi, in preparation.
\end{enumerate}}
 \end{thm}

From now on, we focus on the case when $X/\F_q$ is an elliptic curve $E/\F_q$. By Theorem\,\ref{thm2}(2),  $$\zeta_{E/\F_q,n}(s)=\frac{P_{E/\F_q,n}(T)}{(1-T)(1-QT)}$$ for a certain degree 2 polynomial $P_{E/\F_q,n}(T)$ with rational coefficients. Therefore, there exists a rational number $a_{E/\F_q,n}$ such that
$$P_{E/\F_q,n}(T)=:\a_{E/\F_q,n}\Big(1-a_{E/\F_q,n}T_n+Q_nT_n^2\Big).
$$
Indeed,  the constant term of $P_{E/\F_q,n}(T)$ can be easily seen to coincide with the $\a$-invariant $\a_{E/\F_q,n}(0)$ of $E/\F_q$ from the definition, while the coefficient of $T_n$ in $\frac{1}{\a_{E/\F_q,n}}P_{E/\F_q,n}(T)$ is simply $Q_n$ is a direct consequence of the functional equation, i.e. Theorem\,\ref{thm2}(3).

This rational number $a_{E/\F_q,n}$ introduced above will be called the {\it $a$-invariant in rank $n$ for $E/\F_q$}, since it  clearly is the rank $n$ analogue of the classical $a$-invariant $a_{E/\F_q}:=q+1-\#E(\F_q)$, or equivalently by 
$P_{E/\F_q}(t)=:1-a_{E/\F_q}t+qt^2.$

In the sequel, these $a$-invariants $a_{E/\F_q,n}$ in rank $n$ will be the central players of our studies. Simply put, the existing murmurations for $a_{E/\F_{p_i}}$'s associated to families of elliptic curves
$\E/\Q$'s through their $p_i$-reductions $E/\F_{p_i}$, and the Sate-Tate conjecture for the  $a_{E/\F_{p_i}}$'s associated to these elliptic curves $E/\F_{p_i}$'s, stand equally well for the corresponding  $a$-invariants $a_{E/\F_{p_i},n}$'s in rank $n$.

To state our main results, we make the following preparations:
Let $\E$ be a (regular integral) elliptic curve defined over the field $\Q$ of rationals, and for the $i$-th prime integer $p_i$ $(i\geq 1)$ e.g. $p_1=2, p_2=3,\ldots$, let $E/\F_{p_i}$ denotes the $p_i$-reduction of $\E$. Introduce the rank $n$ average value $f_{r,n}(i)$ by:
\be
f_{r,n}(i):=\frac{1}{\#\cE_r[N_1,N_2]}\times\!\!\!\!\sum_{E\in\cE_r[N_1,N_2]}\bc a^{~}_{E/\F_{p_i},1}&n=1\\
a^{~}_{E/\F_{p_i},2}+p_i-1&n=2\\
\frac{1}{n-1}\cdot\big(a^{~}_{E/\F_{p_i},n}+(n-1)p_i+n-5\big)&n\geq 3\ec
\ee
where $N_1, N_2 \in \Z^+$ satisfying $N_1\le N_2$, and $\cE_r[N_1,N_2]$ denotes the set of elliptic curves over $\Q$ of arithmetic rank $r$ with the conductor in the interval $[N_1, N_2]$.\footnote{Here as in the rank one case,  for each isogeny class of elliptic curves $\E/\Q$, only a single representative elliptic curve is selected in $\cE_r[N_1,N_2]$. 
For examples, as in \cite{MUR1}, from the database of elliptic curves listed in the LMFDB, we  know that 
$$\begin{aligned}
&\#\cE_0[7500,10000]=4238, \ \ \ \#\cE_1[7500,10000]=5194,
\\&\#\cE_0[5000,10000]=8536,\ \ \ \#\cE_2[5000,10000]=1380.
\end{aligned}$$}
\\

In addition, recall that for an elliptic curve $E/\F_q$, by Theorem\,\ref{thm3}, the Riemann hypothesis holds for $\zeta_{E/\F_q,n}$. This is  equivalent to say that 
$$
-1\le\frac{1}{2\sqrt{Q_n}}\cdot a_{E/\F_q,n}\le1.
$$
Since $\cos$-function is strictly decreasing  in the interval $[0,\pi]$, 
accordingly, introduce the rank $n$ argument $\theta_{E/\F_q,n}$ of $E/\F_q$ by
\be
\theta_{E/\F_q,n} := \arccos \Big(\frac{1}{2\sqrt{Q_n}}\cdot a_{E/\F_q,n}\Big)\in [0,\pi].
\ee
At this point, it is pretty tempting to formulate the high rank Sato-Tate distribution using these primitive rank $n$ arguments $\theta_{E/\F_q,n}$. However, one may soon realize that there are several fatal obstacles lying in front.

Indeed, as noticed in \cite{HST}, there are three additional refined structures involved.
Namely, for $n\geq 2$, first,  as the dominant term, the $\theta_{E/\F_q,n}$'s converges to $\frac{\pi}{2}$ when  $q^n$ approaches to infinity. This then yields the so-called {\it small $\delta$-invariant}  in rank $n$ distribution theory of $\E/\F_q$, namely, the Dirac distribution $\delta_{\pi/2}$. Secondly, even after subtracting this accumulating point $\frac{\pi}{2}$, 
the difference appears to be still too tiny to be detected. Hence a suitable huge magnification, namely, $\frac{\sqrt{q^{n-1}}}{n-1}$, should be introduced. But this give rise to the third structure, since
there is still a certain blow-up of level $\frac{1}{2}(\sqrt{q}+\frac{1}{\sqrt{q}})$ for the magnified $\frac{\sqrt{q^{n-1}}}{n-1}\Big(\theta_{E/\F_q,n}-\frac{\pi}{2}\Big)$'s. To control all of them, a new genuine normalization process should be introduced, in order to furmulate our high rank Sato-Tate. With some histories explained in the footnote below, this then finally leads to the current {\it big $\Delta$-invariant}\footnote{In \cite{HST}, a slight different normalization, namely, $\frac{\sqrt{q^{n-1}}}{n-1}(\frac{\pi}{2}-\theta_{E/\F_p,n})+\frac{1}{2}(\sqrt{q}+\frac{1}{\sqrt{q}})$ is introduced. Even, for asymptotic considerations, particularly, for the rank $n$ Sato-Tate conjecture,   it works. But, for smaller $q^n$,  the quantity $\frac{\sqrt{q^{n-1}}}{n-1}(\frac{\pi}{2}-\theta_{E/\F_p,n})+\frac{1}{2}(\sqrt{q}+\frac{1}{\sqrt{q}})$
 may well located outside the closed interval $[-1,1]$.
The idea of the introduction of an additional normalization proves to be very powerful: It leads to its companion for the distributions of high rank zeta zeros of number fields, which itself lately motivated M. Suzuki's work on the rank two zeta zeros of $\Q$: \lq Nearest neighbor spacing distributions for the zeros of the real or imaginary part of the Riemann $\xi$-function on vertical lines', Acta Arith. 170 (2015), 47-65.} $\Delta_{E/\F_q,n}$ {\it in rank $n$ of $E/\F_q$}, defined by:
\be
\Delta_{E/\F_q,n}:=\left\{
\begin{array}{lcl}
\sqrt{q}\cos{\theta_{E/\F_q,2}}+\frac{1}{2}(\sqrt{q}-\frac{1}{\sqrt{q}}) & \mbox{for}& n=2
\\ 
&\\
\frac{\sqrt{q^{n-1}}}{n-1}(\frac{\pi}{2}-\theta_{E/\F_q,n})+\frac{1}{2}(\sqrt{q}+\frac{n-5}{(n-1)\sqrt{q}})  & \mbox{for} & n\ge3 
\end{array}
\right.
\ee

Now we are ready to state the main theorem of this paper:

\begin{thm}\label{thm4} Fix a natural number $n\geq 2$.

\begin{enumerate}
\item [(1)] (Rank $n$ Murmurations) For families of a regular (integral) elliptic curves $\E/\Q$'s, when plotting the points $(i,f_{r,n}(i))$ $i\ge 1$ in the sufficiently large range, the murmuration phenomenon  appear in exactly the same way as the one associated to the $(i,f_{r,1}(i))$'s (of the same families).
\item[(2)] (Rank $n$ Sato-Tate Conjecture) Let $\E/\Q$ be a non CM  elliptic curve. For $\alpha, \beta\in \R$ satisfying $0\le \alpha <\beta \le \pi$, we have
$$
\lim_{N\rightarrow\infty}\frac{\#\{p\le N : p:\ {\rm prime},\ \cos\alpha\ge\Delta_{E/\F_p,n}\ge\cos\beta\}}{\#\{p\le N: p:\ {\rm prime}\}}=\frac{2}{\pi}\int_\alpha^\beta\sin^2\theta d\theta.
$$
\end{enumerate}
\end{thm}

To entertaining  the reader, we next provide  some concrete examples on the rank $n$ murmurations with fixed ranges on the  conductors, the arithmetic rank $r$ of elliptic curves $\E/\Q$ and the geometric rank $n$.

\begin{figure}[H]
    \centering
    \includegraphics[width=10.0cm]{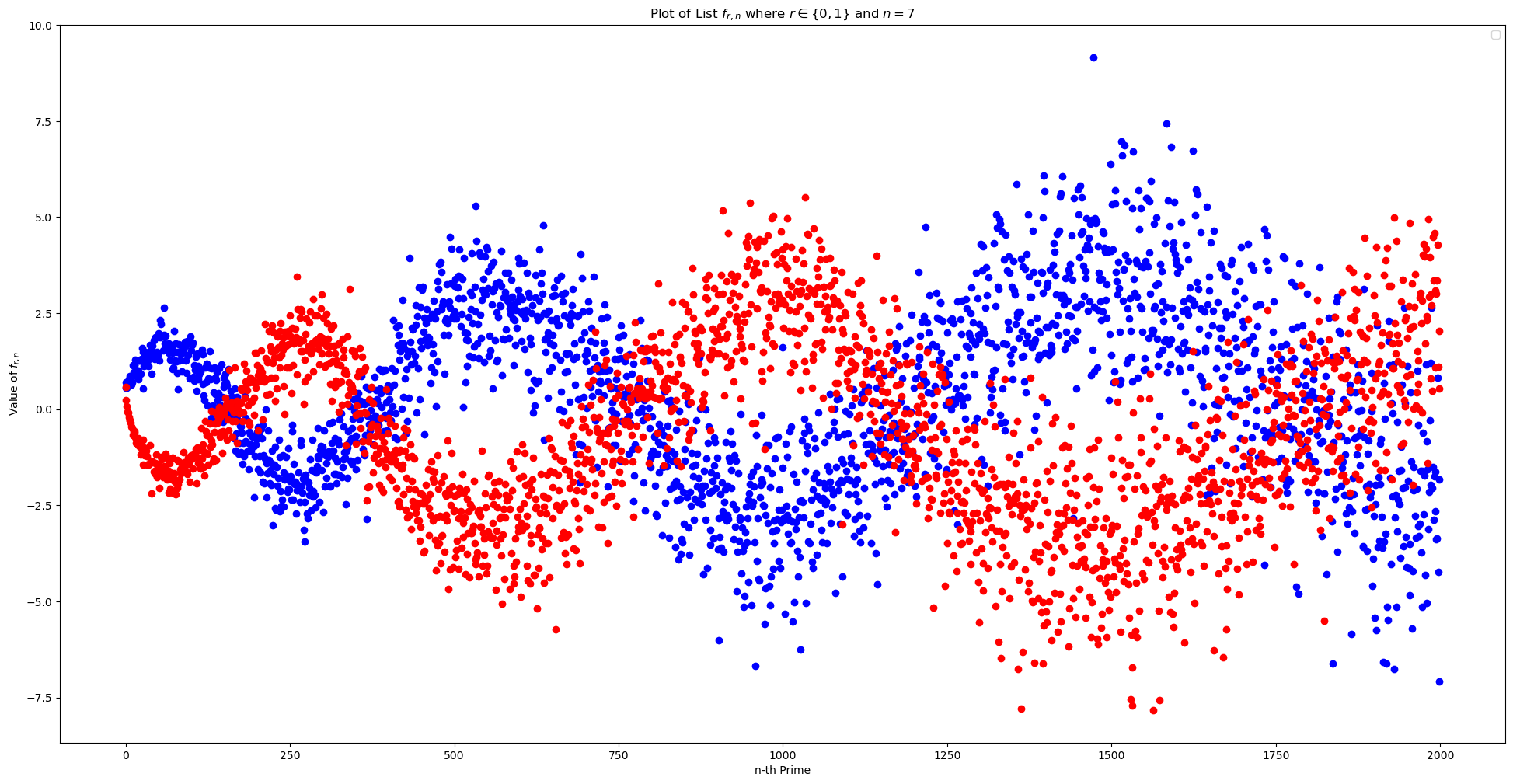}
    \caption{Plot of $f_{r,n}(i)$ where $r\in{0,1}$ and $n=7$, for elliptic curves with conductor in $[7500,10000]$. $f_{0,n}(i)$ is in blue and $f_{1,n}(i)$ is in red.}
    \label{hst}
\end{figure}

\begin{figure}[H]
    \centering
    \includegraphics[width=10.0cm]{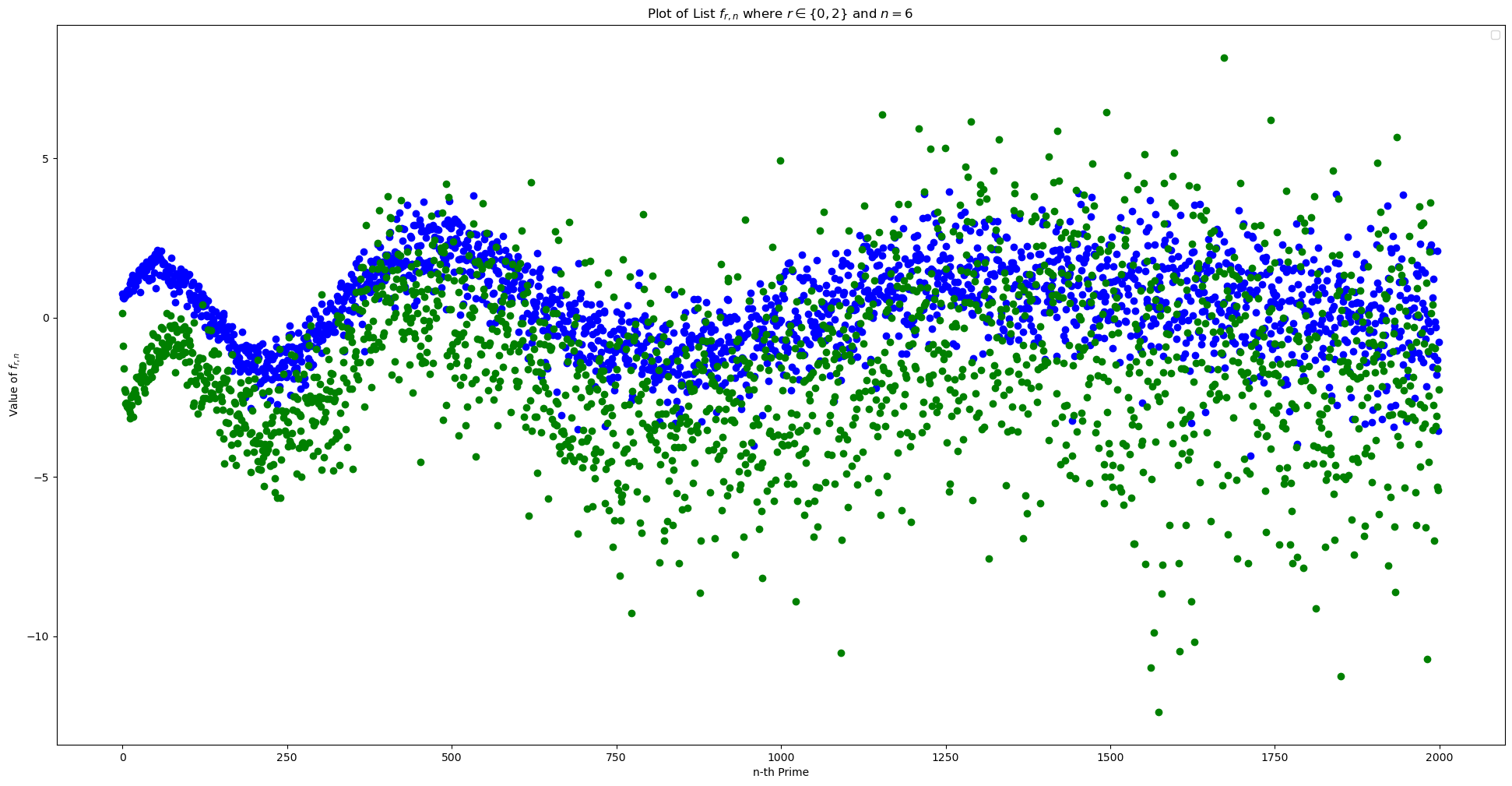}
    \caption{Plot of $f_{r,n}(i)$ where $r\in{0,2}$ and $n=6$, for elliptic curves with conductor in $[5000,10000]$. $f_{0,n}(i)$ is in blue and $f_{2,n}(i)$ is in green.}
    \label{hst}
\end{figure}

As for the rank $n$ Sato-Tate distribution for the  big $\Delta$-invariants $\Delta_{E/\F_q,n}$ in rank $n$ of elliptic curve $\E/\Q$, we provide the following illustrative examples.

\begin{figure}[H]
    \centering
    \includegraphics[width=8.0cm]{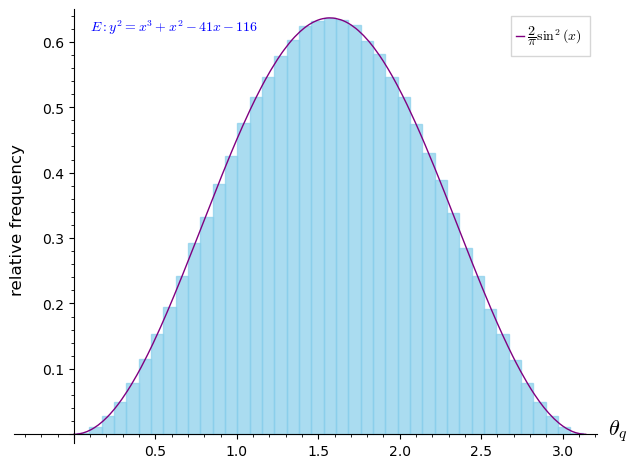}
    \caption{Sato-Tate distribution of rank 3 zeta function $\zeta_{E/\F_q,3}(s)$ over elliptic curve $\E/\Q: y^2 = x^3 + x^2 - 41x - 116$ and $q \le N =10,000,000$.}
    \label{hst}
\end{figure}

\begin{figure}[H]
    \centering
    \includegraphics[width=7.0cm]{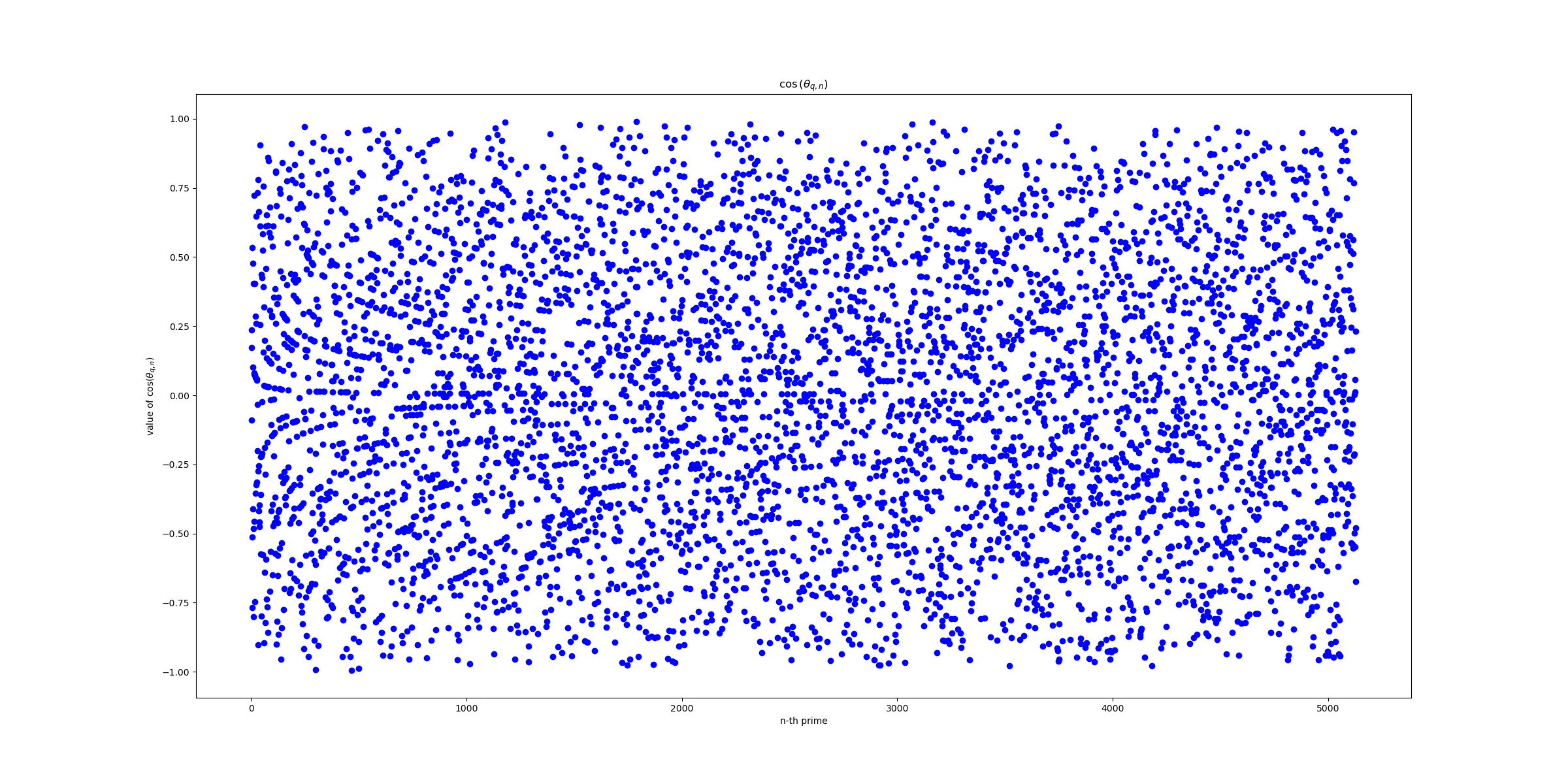}
    \caption{Plot of $\Delta_{E/\F_q,n}$ over elliptic curve $E: y^2 = x^3 + x^2 - 41x - 116$ and $q\le N=50,000$ when $n=5$.}
    \label{cos5w}
\end{figure}

\begin{figure}[H]
    \centering
    \includegraphics[width=8.0cm]{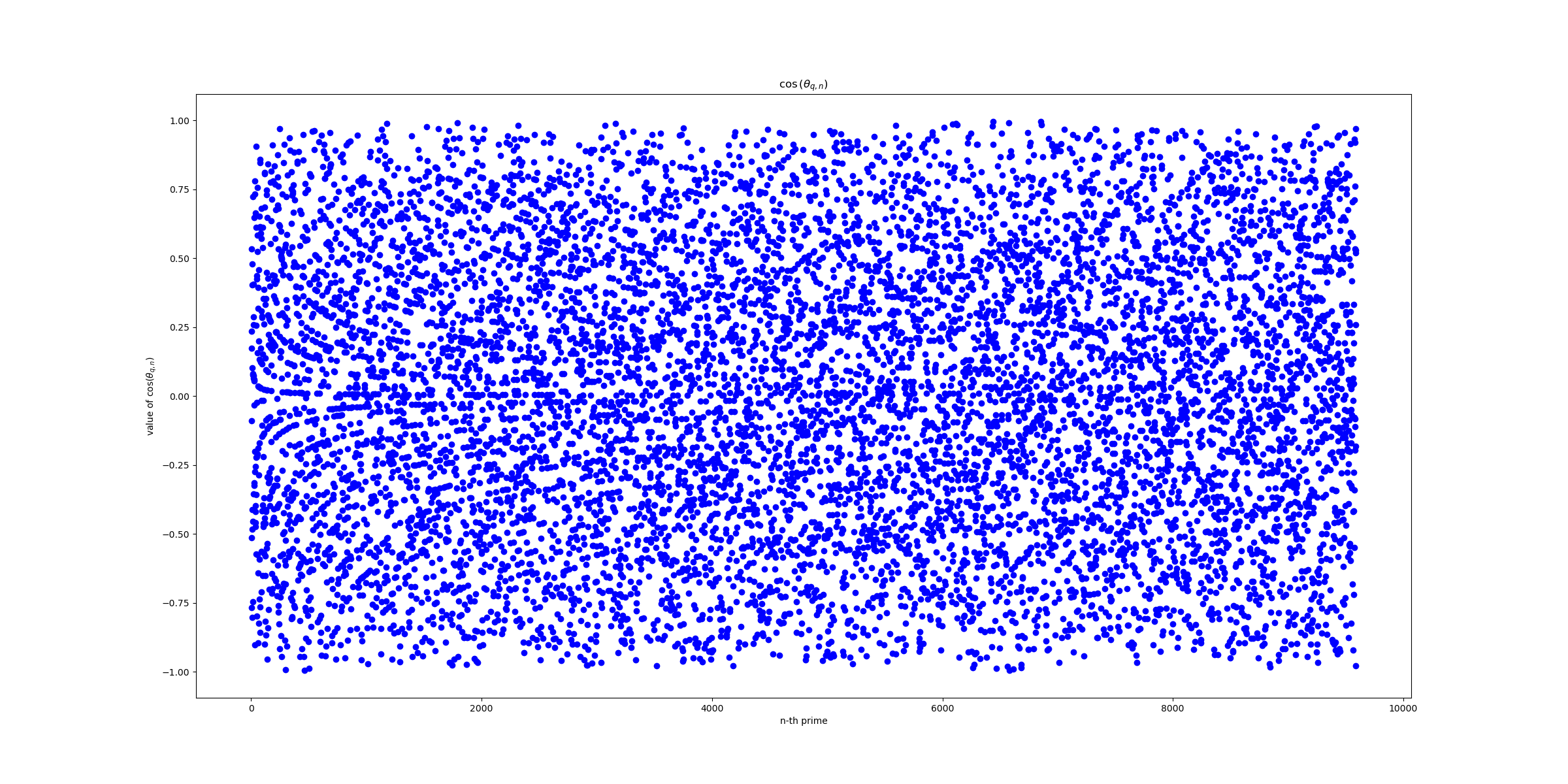}
    \caption{Plot of $\Delta_{E/\F_q,n}$ over elliptic curve $E: y^2 = x^3 + x^2 - 41x - 116$ and $q\le N=100,000$ when $n=5$.}
    \label{cos10w}
\end{figure}

\begin{figure}[H]
    \centering
    \includegraphics[width=8.0cm]{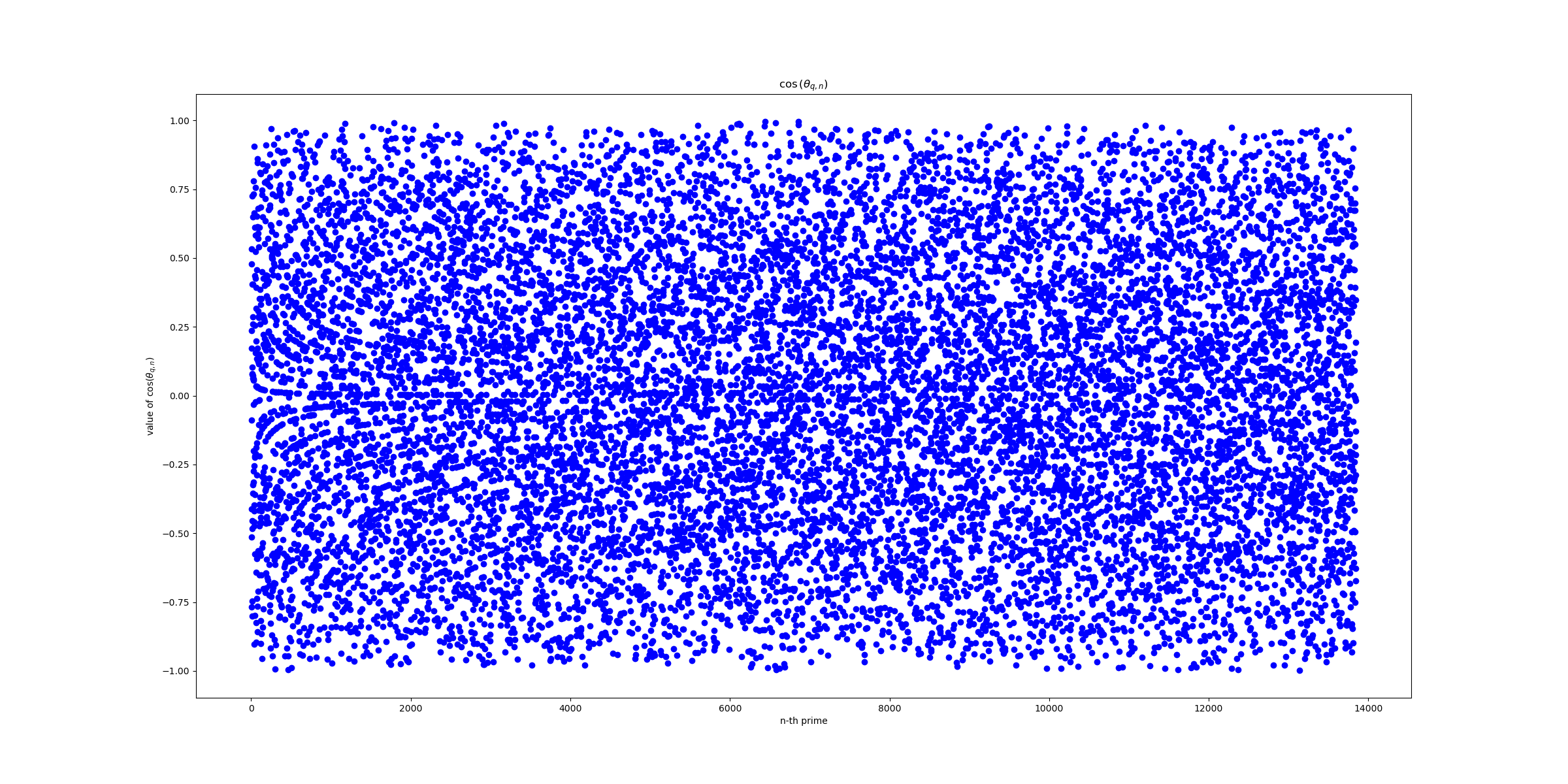}
    \caption{Plot of $\Delta_{E/\F_q,n}$ over elliptic curve $E: y^2 = x^3 + x^2 - 41x - 116$ and $q\le N=150,000$ when $n=5$.}
    \label{cos15w}
\end{figure}

\begin{figure}[H]
    \centering
    \includegraphics[width=8.0cm]{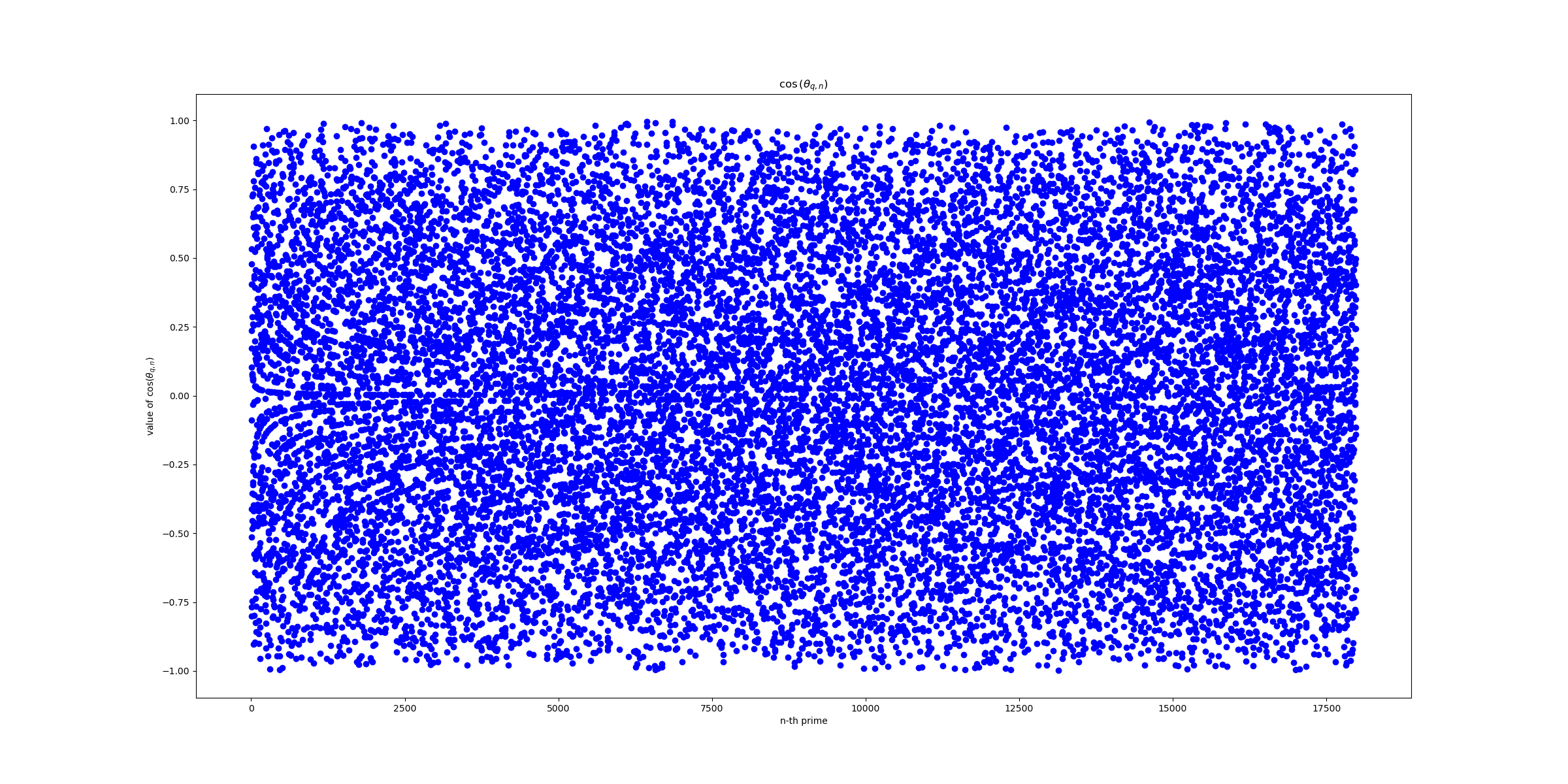}
    \caption{Plot of $\Delta_{E/\F_q,n}$ over elliptic curve $E: y^2 = x^3 + x^2 - 41x - 116$ and $q\le N=200,000$ when $n=5$.}
    \label{cos20w}
\end{figure}

To end this section, we point out that there are clear geometric flows from up/lower left to the middle right  appeared in the final four pictures above. It would be very interesting to explore the hidden mathematical structures, if any, behind.

\section{Proof of Main Theorem}

From now on, we will focus on the situation when $X/\F_q$ is an (integral  projective regular) elliptic curve $E/\F_q$. Fix a natural number $n$ unless otherwise is stated. 

Since for a semi-stable vector bundle $V$ of rank $n$ on $E/\F_q$, if its degree $d=nm$ is $\gneq 0$, or the same $\geq n$, we have the well-known vanishing theorem $h^1(E,V)=0$. In particular, for such a $V$, $h^0(E,V)=d$, by the Riemann-Roch theorem. Consequently, 
\be\label{eq17}
\alpha_{E/\F_q,n}(d)=(q^d-1)\beta_{E/\F_q,n}(d), \quad \beta_{E/\F_q,n}(mn)=\beta_{E/\F_q,n}(0), \forall d=mn,\ m\in \N.
\ee
Here in the second relation, we have used the fact that each rank $n$ semi-stable vector bundle on $E/\F_q$ of  degree $nm$ is obtained  from a uniquely determind rank $n$ semi-stable vector bundle on $E$ of degree $0$ by tensoring an $m$-th multiple of a certain line bundle on $E/\F_q$ of  degree 1, whose existence is guaranteed by a well-known result of Artin (See e.g. \cite{Sil}). Accordingly, by the well-known Duality, a standard procedure in zeta function discussion, \eqref{eq8} and \eqref{eq17}, implies that
\be\label{eq18}
\zeta_{E/\F_q,n}(s)=\alpha_{E/\F_q,n}(0)+\beta_{E/\F_q,n}(0)\cdot \frac{(Q_n-1)T_n}{(1-T_n)(1-Q_nT_n)}.
\ee

With $E$, $\F_q$ and $n$ fixed for the time being, we simply write
$$
\alpha_n:=\alpha_{E/\F_q,n}(0),\quad \beta_n:=\beta_{E/\F_q,n}(0)\qquad{\rm (and}\qquad Q=Q_n,\quad T=T_n,)
$$
so that the rank $n$ zeta function for elliptic curve $E/\F_q$ and its associated polynomial $P_{E,n}(T)=P_{E/\F_q,n}(T_n)$ can be simply expressed as:
\be
\zeta_{E,n}(s)=\alpha_n+\beta_n\cdot \frac{(Q-1)T}{(1-T)(1-QT)}\qqan P_{E,n}(T)=\alpha_n\Big(1-a_{E,n}T+QT^2\Big).
\ee
where $$a_{E,n}:=a_{E/\F_q,n}:=(Q_n+1)-(Q_n-1)\frac{\beta_{E/\F_q,n}(0)}{\alpha_{E/\F_q,n}(0)}=(Q+1)-(Q-1)\frac{\beta_n}{\alpha_n},$$ as to be seen by a direct calculation from \eqref{eq18}.

The $\alpha$- and $\beta$- invariants have been intensively studied with a long history and a rich theory (see e.g. \cite{EC}, \cite{SLn}, \cite{HST} and \cite{SGHR},\cite{MR}). In fact, for a general curve $X/\F_q$, $\beta_{X/\F_q,n}$ was first introduced by 
Harder-Narasimhan in \cite{HN}.  Here for our limited purpose, we select the follows.
\begin{thm}[$\alpha$- and $\beta$- invariants, \cite{EC}]
For an elliptic curve $E/\F_q$,  we have:
\begin{enumerate}
\item[(1)] (Counting Miracle) For all $n\ge0$,  
\be
\alpha_{n+1}=\beta_n.
\ee
\item[(2)] (Structural Recursion Formula) For  all $n\ge1$, the $\beta_n$'s satisfy a two-step recursion relation
\be
(q^n-1)\beta_n=(q^n+q^{n-1}-a_{E/\F_q,1})\beta_{n-1}-(q^{n-1}-q)\beta_{n-2},
\ee
with initial values 
$\beta_0=1$ and $\beta_{-1}=0$.
\end{enumerate}
\end{thm}
\begin{rem} \normalfont{The counting miracle was first conjectured by the senior author for elliptic curves based on some detailed calculations in lower ranks, with the help of the so-called Atiyah bundles. As indicated above, the so-called Counting Miracle Conjectured was first proved in (the first version of) \cite{EC} by Weng-Zagier for elliptic curves $E/\F_q$. Later, this was generalized for general (integral regular) projective curves $X/\F_q$ of genus $g$
by Sugahara (\cite{SGHR}) and  independently  by Mozgovoy-Reineke (\cite{MR}) in the form:}
$$
\alpha_{X/\F_q,n+1}(0)=q^{n(g-1)}\beta_{X/\F_q,n}(0).
$$
\end{rem}

Since the counting miracle offers us an intrinsic relation between the $\alpha$- and $\beta$-invariants, all the higher rank zeta functions for elliptic curves (over finite fields) are totally determined by the sequence $\{\b_n\}_{n\geq -1}$ of the beta invariants. That is to say, for $n\geq 1$,
 the rank $n$ zeta function of elliptic curve $E/\F_q$ is given by
\bea
\zeta_{E/\F_q,n}(s)
=&\beta_{n-1}+\beta_n\cdot \frac{(Q-1)T}{(1-T)(1-QT)}.
\eea
Consequently, the polynomial $P_{E/\F_q,n}(T)$ simply takes the form
\be\label{eq23}
\frac{1}{\b_{n-1}}P_{E/\F_q,n}(T)=1-\left((Q+1)-(Q-1)\frac{\beta_n}{\beta_{n-1}}\right)T+QT^2,
\ee
and the associated  $a$-invariant in rank $n$ is given by
$$a_{E/\F_q,n}=(Q+1)-(Q-1)\frac{\beta_n}{\beta_{n-1}}.$$

With this point settled, next we turn to the simple 2-step recursion formula on the $\b_n$'s, from which the $\b_n$ in rank $n$ can be determined  completely from its initial values $\b_{-1}=0$ and $\b_0=1$ in terms of $a_{E/\F_q,1}$ and $q, n$.

\begin{ex} \normalfont{When $n=1$, we have
$$\begin{aligned}
(q^1-1)\beta_{E/\F_q,1}
=&(q^1+q^{1-1}-a_{E/\F_q,1})\beta_{1-1}-(q^{1-1}-q)\beta_{1-2}\\
=&q+1-a_{E/\F_q,1}=\#E(\F_q).\end{aligned}
$$
Accordingly, the rank one zeta function of $E/\F_q$ becomes simply 
$$
\zeta_{E,1}(s)=\beta_0+\beta_{E/\F_q,1}\cdot \frac{(q^1-1)t^1}{(1-t^1)(1-q^1t^1)}
=\frac{1-a_{E/\F_q,1}t+qt^2}{(1-t)(1-qt)}$$
which is nothing but the classical Hasse-Weil zeta function $\zeta_{E/\F_q}(s)$ of $E/\F_q$.

Similarly, when $n=2$, we have 
$$\begin{aligned}
(q^2-1)\beta_2=&(q^2+q^{2-1}-a_{E/\F_q,1})\beta_{2-1}-(q^{2-1}-q)\beta_{2-2}\\
=&\frac{(q^2+q-a_{E/\F_q,1})(q+1-a_{E/\F_q,1})}{q-1}.\end{aligned}
$$
Accordingly, the rank two zeta function of $E/\F_q$  becomes
$$\begin{aligned}
\zeta_{E,2}(s)
=&\beta_{E/\F_q,1}+\beta_2\cdot \frac{(q^2-1)t^2}{(1-t^2)(1-q^2t^2)}\\
=&\frac{q+1-a_{E/\F_q,1}}{q-1}\times\frac{1-(a_{E/\F_q,1}-q+1)T+QT^2}{(1-T)(1-QT)}
\end{aligned}$$

Obviously, $\alpha_2=(q+1-a_{E/\F_q,1})/(q-1)=\b_1$ is a constant depending merely on the elliptic curve $E/\F_q$ and, in particular,}
$$a_{E,1}=a_{E/\F_q,1}=q+1-\#E(\F_q)\qqan a_{E,2}=a_{E/\F_q,1}-q+1.\hskip 1.0cm\square$$
\end{ex}~\\

Back to the general discussion, 
by the recursion relation, namely Theorem 4(2), we have, for $n>0$,
\be\label{eq17*}
\frac{\beta_n}{\beta_{n-1}}=\frac{q^n+q^{n-1}-a_{E/\F_q,1}}{q^n-1}-\frac{q^{n-1}-q}{q^n-1}\times\frac{\beta_{n-2}}{\beta_{n-1}}.
\ee
To simplify our notations further, for our own convenience, set now 
\be
\gamma_n:=\frac{\beta_n}{\beta_{n-1}}
,\quad
A_n:=\frac{q^n+q^{n-1}-a_{E/\F_q,1}}{q^n-1},
\qan
B_n:=\frac{q^{n-1}-q}{q^n-1}, 
\ee
we have, from \eqref{eq17*} 
\be
\gamma_n=A_n-B_n\gamma_{n-1}^{-1}.
\ee
Consequently, by definition,
\bea\label{eq20}
a_{E/\F_q,n}&=(q^n+1)-(q^n-1)\gamma_n
\\
&=(q^n+1)-(q^n-1)\Big(A_n-B_n\times\frac{q^{n-1}-1}{q^{n-1}+1-a_{E/\F_q,n-1}}\Big).
\eea
where in the last equality we have used the level $(n-1)$ definition
\be
a_{E/\F_q,n-1}=(q^{n-1}+1)-(q^{n-1}-1)\gamma_{n-1}.
\ee

Now we are ready to introduce our main technical result to prove Theorem\,\ref{thm4} on the murmurations and the Sato-Tate conjecture for rank $n$ zeta functions of elliptic curves:
\begin{thm}\label{thm6}[Asymptotic behavior of $a_{E/\F_q,n}$]
We have
\be
a_{E/\F_q,1}=a_{E/\F_q}, \qquad a_{E/\F_q,2}=1+a_{E/\F_q,1}-q,
\ee
and 
\be
a_{E/\F_q,n}=(5-n)+(n-1)a_{E/\F_q,1}-(n-1)q+O\Big(\frac{1}{\sqrt{q}}\Big)\quad (n\geq 3)
\ee
In particular, for $n\geq 3$
\be\label{eq24}
a_{E/\F_q,n}\sim(5-n)+(n-1)a_{E/\F_q,1}-(n-1)q\qquad (q\rightarrow\infty).
\ee
\end{thm}
From this theorem, up to a well-understood  normalization if necessary, the asymptotic behavior of $a$-invariant $a_{E/\F_q,n}$ in rank $n$ is reduced to that of $a_{E/\F_q,1}$. This indicates that the murmuration and the Sato-Tate conjecture for $a_{E/\F_q}$ associated to the classical Hasse-Weil zetas of elliptic curves work in exactly the same way for all high rank zetas as well.
\bp
To prove this theorem, we use an induction on $n$, based on \eqref{eq20} and the standard Hasse-Weil bound for elliptic curves, i.e. $a_{E/\F_q}=O(\sqrt{q})$.

Indeed, by \eqref{eq20}, we have
\bea\label{eq25}
a_{E/\F_q,n}=&(q^n+1)-(q^n-1)\frac{q^n+q^{n-1}-a_{E/\F_q,1}}{q^n-1}
\\&\qquad+(q^n-1)\frac{q^{n-1}-q}{q^n-1}\times\frac{q^{n-1}-1}{q^{n-1}+1-a_{E/\F_q,n-1}}
\\
&=1+a_{E/\F_q,1}+\frac{-q^n-2q^{n-1}+a_{E/\F_q,n-1}q^{n-1}+q}{q^{n-1}+1-a_{E/\F_q,n-1}}.
\eea

We start with $n=3$. From the recursion relation \eqref{eq25} above, we have 
\bea\label{eq26}
a_{E/\F_q,3}
&=1+a_{E/\F_q,1}+\frac{-q^3-2q^2+a_{E/\F_q,2}q^2+q}{q^2+1-a_{E/\F_q,2}}\\
&=1+a_{E/\F_q,1}+\frac{-2q^3+(a_{E/\F_q,1}-1)q^2+q}{q^2+q-a_{E/\F_q,1}}
\\&\hskip 4.0cm({\rm using}\ a_{E/\F_q,2}=1+a_{E/\F_q,1}-q)\\
&=1+a_{E/\F_q,1}+\frac{-2q(q^2+q-a_{E/\F_q,1}-1)}{q^2+q-a_{E/\F_q,1}}\\
&\qquad\quad\qquad\quad+\frac{(2+a_{E/\F_q,1}-1)q^2+O(a_{E/\F_q,1}\cdot q)}{q^2+q-a_{E/\F_q,1}}
\\
&\hskip 4.0cm({\rm make\ an\ oriented\ seperation})\\
&=1+a_{E/\F_q,1}-2q+1+a_{E/\F_q,1}+\frac{O(a_{E/\F_q,1}\cdot q)}{q^2+q-a_{E/\F_q,1}}
\\&=2+2a_{E/\F_q,1}-2q+O\Big(\frac{1}{\sqrt{q}}\Big)
\eea
as wanted. Here in the last three equalities, we have used the Hasse bound $a_{E/\F_q,1}=O\Big(\sqrt q\Big)$.

Assume now inductively that the assertion holds for level $n\geq 3$, that is,  
 $$a_{E/\F_q,n}=(5-n)+(n-1)a_{E/\F_q,1}-(n-1)q+O\Big(\frac{1}{\sqrt{q}}\Big)\qquad n\geq 3.$$ 
 To understand $a_{E/\F_q,n+1}$, using the recursion relation \eqref{eq25}, we have, as $q\rightarrow\infty$,
$$\begin{aligned}
&a_{E/\F_q,n+1}\\
&=1-q^n+a_{E/\F_q,1}+\frac{(q^n-q)(q^n-1)}{q^n+1-a_{E/\F_q,n}}
\\&=1-q^n+a_{E/\F_q,1}+\frac{(q^n-q)(q^n-1)}{q^n+(n-1)q-(n-1)a_{E/\F_q,1}-(5-n)+O\Big(\frac{1}{\sqrt{q}}\Big)}
\\&=1+a_{E/\F_q,1}+\frac{-nq^{n+1}+((4-n)+(n-1)a_{E/\F_q,1}+O\Big(\frac{1}{\sqrt{q}}\Big))q^n}{q^n+(n-1)q-(n-1)a_{E/\F_q,1}-(5-n)+O\Big(\frac{1}{\sqrt{q}}\Big)}
\\&=(5-(n+1))+na_{E/\F_q,1}-nq+O\Big(\frac{1}{\sqrt{q}}\Big).
\end{aligned}$$
This is exact what we are looking after.
\ep

Before leaving here, we point out that, for a fixed $n$, $a_{E/\F_q,n}$ is always of the order $O(q)$. This means that in \eqref{eq25}, the denominator is of order $q^{n-1}+O(q)$. This also explains why in the final form of \eqref{eq26}, the constant term is 2 instead of $-1$ appeared in the recursion relation \eqref{eq25}. ~\\

\noindent
{\it Proof of Theorem\,\ref{thm4}(1) for murmurations of rank $n$ zeta functions}.
Murmurations were first discovered  by He-Lee-Oliver-Pozdnyako in \cite{MUR1}, using plots of average values $f_{r,1}(i)$ of the geometric rank one coefficients $a_{E/\F_{p_i},1}$ of the $p_i$-reductions $E/\F_{p_i}$ for families of elliptic curves $\E/\Q$ of the given arithmetic ranks $r$ and with conductors in a fixed range. Based on these average values, the elliptic curves are grouped  beautifully according to their arithmetic ranks $r$. Some data-scientific experimentations and machine learnings have been performed in \cite{MUR1}. Next, we establish the murmurations for high rank zeta functions of elliptic curves similarly based on Theorem\,\ref{thm6}. 

Let $\E/\Q$ be a given elliptic curve, and let $p_i$ be the $i$-th prime integer, e.g. $p_1=2, p_2=3,\ldots$. Then for a fixed (geometric rank) $n\ge2$, we examine the rank $n$ zeta function 
\be
\zeta_{E/\F_{p_i},n}(s)=\frac{1}{\alpha_n}\times\frac{1-a_{E/p_i,n}T+p_i^nT^2}{(1-T)(1-p_i^nT)},
\ee
where $E/\F_p$ denotes the $p$-reduction of $\E/\Q$.
We have introduced the following average values  ($i=1,2,\ldots$):
\be
f_{r,n}(i):=\frac{1}{\#\cE_r[N_1,N_2]}\times\!\!\!\!\sum_{E\in\cE_r[N_1,N_2]}\bc a^{~}_{E/\F_{p_i},1}&n=1\\
a^{~}_{E/\F_{p_i},2}+p_i-1&n=2\\
\frac{1}{n-1}\cdot\big(a^{~}_{E/\F_{p_i},n}+(n-1)p_i+n-5\big)&n\geq 3\ec
\ee
where $N_1, N_2 \in \Z^+$ satisfying $N_1\le N_2$, and $\cE_r[N_1,N_2]$ denotes the set of elliptic curves over $\Q$ of rank $r$ whose conductors belong to the interval $[N_1, N_2]$. 

Note that by Theorem\,\ref{thm6}, particularly, \eqref{eq24}, we have, for $q=p^k$,
$$a_{E/\F_q,n}=\bc a_{E/\F_q,1}&n=1\\
1+a_{E/\F_q,1}-q&n=2\\
\sim (5-n)+(n-1)a_{E/\F_q,1}-(n-1)q\qquad (q\rightarrow\infty)&n\geq 3\ec.$$
Hence our result comes as a direct consequence of the following main result of \cite{MUR1} for the murmurations for elliptic curves in terms  of the $a_{E/\F_{p_i}}$'s.

\begin{obs}[Murmurations for elliptic curves, \cite{MUR1}]
When plotting the points $(i,f_{r,1}(i))$ for $i$ from $1$ to a large integer for families of elliptic curves $\E/\Q$, murmurations do appear according to the arithmetic ranks $r$.
\end{obs}

We end this discussions with the following comments. Murmurations of high rank zeta function have been observed to form the same pattern as these using rank one Hasse-Weil zeta functions for elliptic curves. This indicates that many techniques and skills in the theory of  Hasse-Weil zeta functions can also be equally applied to the high rank zeta functions. For instance, as in \cite{MUR1}, the Principal Component Analysis(PCA) has been preformed  to a balanced set of $36,000$ randomly selected elliptic curves of $r\le2$, with $N_1=1\times 10^4$ and $N_4=4\times 10^4$. Moreover, for the case $r\in\{0,1\}$, curves of the form
\be
y=Ax^\alpha\sin(Bx^\beta)
\ee
can be used to fit the data on the average of $a_{E/\F_{p_i},1}$. In addition, some more complicated experiments about murmurations are introduced in \cite{MUR2}. Sincerely hope that our high rank zeta functions will provide various mathematics communities some new test grounds to study mathematical structures for various curves, among others.\hskip 8.60cm$\square$~\\

\noindent
{\it Proof of Theorem\,\ref{thm4}(2) on the rank $n$ Sato-Tate Conjecture}.
Our proof is based on Theorem\,\ref{thm6} above and the Barnet(-)Lamb-Geraghty-Harris-Taylor's Theorem on the original Sato-Tate conjecture. For this reason, it is more convenient for us to recall
the classical Sato-Tate conjecture for Hasse-Weil zeta function of elliptic curves $E/\F_{p_i}$'s.

Let $\E/\Q$ be an (integral regular projective) elliptic curve and denote its $p$ reduction by $E/\F_p$. Then, by
the Riemann hypothesis for the classical Hasse-Weil zeta function for the  elliptic curve $E/\F_p$, 
\be
-1<\frac{a^{~}_{E/\F_p,1}}{2\sqrt{p}}=\frac{1+p-\#E(\F_p)}{2\sqrt{p}}<1,
\ee
and thus we may introduce the associated argument 
$\theta_{E/\F_p,1}\in [0,\pi]$ by 
\be
\cos\theta_{E/\F_p,1}=\frac{a_{E/\F_p,1}}{2\sqrt{p}}.
\ee

\begin{thm}[Sato-Tate Conjecture \cite{BLGHT}, \cite{CHT}, \cite{HT}, \cite{T}, see also \cite{ST1} and \cite{ST2} for backgrounds]\label{thm8} Let $\E/\Q$ be a fixed non CM elliptic curve. 
Then for any real $\a,\b$ satisfying $0\le \alpha <\beta \le \pi$, the following holds:
\be
\lim_{N\rightarrow\infty}\frac{\#\{p\le N : p \ {\rm prime}, \alpha\le\theta_{E/\F_p,1}\le\beta\}}{\#\{p\le N: p \ {\rm prime}\}}=\frac{2}{\pi}\int_\alpha^\beta\sin^2\theta d\theta.
\ee
\end{thm}

Next, we consider its higher rank analogue. By a result of  Weng-Zagier \cite{EC}, e.g. Theorem\,\ref{thm3}, the rank $n$ zeta function
 $\zeta_{E/\F_q,n}(s)$ of an elliptic curve $E/\F_q$ satisfies the Riemann hypothesis as well.
 This is equivalent to 
\be
-1\le\frac{1}{2\sqrt{Q_n}}\cdot a_{E/\F_q,n}\le1
\ee
Accordingly,  we may introduce the associated argument 
$\theta_{E/\F_q,n}\in [0,\pi]$ via
$$\cos\theta_{E/\F_q,n}=\frac{1}{2\sqrt{Q_n}}\cdot a_{E/\F_q,n}\quad{\rm or\ equivalently,}\qquad 
\theta_{E/\F_q,n} := \arccos \frac{a_{E/\F_q,n}}{2\sqrt{Q}}.$$

It is quite tempting to formulate the rank $n$ Sato-Tate conjecture as in rank 1 case by directly using
$\theta_{E/\F_p,n}$. The fact is that the structures here are rather subtle. Indeed, as exposed in an earlier work by the senior author (see e.g. \cite{HST}), based on a rough asymptotic relation  on the $\gamma$-invariant of elliptic curve $E/\F_q$ in (\cite{EC}), namely, 
\be
\gamma_n:=\frac{\beta_n}{\beta_{n-1}}=1+\frac{(n-1)(q-a_{E/\F_q,1}+1)-c(q)}{q^n}+O\Big(\frac{n^2}{q^{2n-2}}\Big)
\ee
as $q^n\to \infty$, where $c(q)=2+3(a-2)/q+...$ is independent of $n$. In particular,
 we have the following earlier result, which can now be verified as a direct consequence of Theorem\,\ref{thm4}(2): 
 
\begin{thm*}[First Approximation to Rank $n(\geq 3)$ Sato-Tate Distribution, see \cite{HST}]
For a non CM elliptic curve $\E/\Q$, with respect to real numbers $\a,\b$ satisfying $0\le \alpha <\beta \le \pi$, we have:
\bea
&\lim_{N\rightarrow\infty}
\frac{\#\{p\le N : p\ {\rm prime},\ \cos\alpha\,\dot\ge\frac{\sqrt{p^{n-1}}}{n-1}\Big(\frac{\pi}{2}-\theta_{E/\F_p,n}\Big)\!+\!\frac{1}{2}\Big(\sqrt{p}+\frac{1}{\sqrt{p}}\Big)\dot\ge\,\cos\beta\}}{\#\{p\le N: p\ {\rm prime}\}}\\
&=\frac{2}{\pi} \int_\alpha^\beta {\sin}^2 \theta d\theta.
\eea
\end{thm*}

Our proof of Theorem\,\ref{thm4}(2) below is certainly independent of his result. But we recall it here for the purpose to   explain the complicated refined structures on the distributions of $\theta_{E/\F_p,n}$ involved. Indeed, this result implies that there are three levels of such refined structures on the distributions of $\theta_{E/\F_p,n}$. Namely, the first on the small Dirac delta distribution  $\delta_{\pi/2}$, which means that the $\theta_{E/\F_p,n}$'s
have an accumulate point at $\frac{\pi}{2}$, and  even after this limit point is subtracted, to witness the secondary structure, an essential magnified  scale factor $\frac{\sqrt{p^{n-1}}}{n-1}$ should be introduced. However, with such a huge magnification introduced, a blowing-up has been unavoidably
introduced at the same time. Hence thirdly, we should control this blow-up by subtracting  $\frac{1}{2}\Big(\sqrt{p}+\frac{1}{\sqrt{p}}\Big)$. With all these normalizations introduced, then the rank $n$ Sato-Tate distribution can be introduced and proved to perform in  exactly  the same as the classical one for the rank one abelian theory.

Even this earlier version works nicely in asymptotic sense which is enough for the rank $n$ Sato-Tate, 
there is a small defect in this discussion. Namely, there may well exist some smaller $p$'s such that the middle term appeared above, i.e. $\frac{\sqrt{p^{n-1}}}{n-1}(\frac{\pi}{2}-\theta_{E/\F_p,n})+\frac{1}{2}(\sqrt{p}+\frac{1}{\sqrt{p}})$ may well lie outside the closed interval $[-1,1]$. This is the reason why a dot above the inequalities on both sides is introduced, in order to ignore all the contributions that are out of the range $[-1,1]$.

This defect has been remedy in the present work: while asymptotically both are the same, but the current one stands nicely. Indeed, the rank $n$ Sato-Tate conjecture, as stated  in Theorem\,\ref{thm4}(2), now reads as follows: 
Let  $\E/\Q$ be a non CM elliptic curve. Then,  for any two real numbers $\a,\b$ satisfying  $0\le \alpha <\beta \le \pi$, we have:
\be
\lim_{N\rightarrow\infty}\frac{\#\{p\le N : p\ {\rm prime},\ \cos\alpha\ge\Delta_{E/\F_p,n}\ge\cos\beta\}}{\#\{p\le N:  p\ {\rm prime}\}}=\frac{2}{\pi}\int_\alpha^\beta\sin^2\theta d\theta,
\ee
where $\Delta_{E/\F_p,n}$ is defined as:
\be
\Delta_{E/\F_p,n}:=\left\{
\begin{array}{lcl}
\sqrt{p}\cos{\theta_{E/\F_p,2}}+\frac{1}{2}\Big(\sqrt{p}-\frac{1}{\sqrt{p}}\Big) & \mbox{for}& n=2
\\ \frac{\sqrt{p^{n-1}}}{n-1}\Big(\frac{\pi}{2}-\theta_{E/\F_p,n}\Big)+\frac{1}{2}\Big(\sqrt{p}+\frac{n-5}{(n-1)\sqrt{p}}\Big)  & \mbox{for} & n\ge3 
\end{array}
\right.
\ee

In the sequel,
we only consider when $n\ge 3$ since the case for $n=2$ is trivial.
We start with the following
\begin{lem}[Small $\delta$ distribution] \label{lem 7} With the same notation as above, the $\theta_{E/\F_p,n}$'s converge to $\frac{\pi}{2}$ as $p\to\infty$
\end{lem}
\bp
From the asymptotic behaviors of $a_{E/\F_p,n}$'s in Theorem\,\ref{thm6} and the classical Hasse theorem for Hasse-Weil zeta function of $\E/\F_p$, $a_{E/\F_p,n}=O(p)$. This, together with the rank $n$ Riemann hypothesis for $E/\F_q$, namely, Theorem\,\ref{thm3} implies that $\theta_{E/\F_p,n}$ has a limit at $\frac{\pi}{2}$ as $p\to\infty$, which yields  the small Dirac delta distribution $\delta_{\pi/2}$. Similarly,  we have $\a_{E/\F_p,n}/2\sqrt{p^n}$ goes to $0$ when $p$ approaches to $\infty$.
\ep

By definition,  (with the lemma,) we have
\bea
\Delta_{E/\F_p,n}
&=\frac{\sqrt{p^{n-1}}}{n-1}\Big(\frac{\pi}{2}-\arccos \frac{a_{E/\F_p,n}}{2\sqrt{p^n}}\Big)+\frac{1}{2}\Big(\sqrt{p}+\frac{n-5}{(n-1)\sqrt{p}}\Big)
\\&=\frac{\sqrt{p^{n-1}}}{n-1}\Big(\arcsin \frac{a_{E/\F_p,n}}{2\sqrt{p^n}}\Big)+\frac{1}{2}\Big(\sqrt{p}+\frac{n-5}{(n-1)\sqrt{p}}\Big).
\eea

Consequently, for sufficiently large $p$, from the lemma above, we have  
$$
\Delta_{E/\F_p,n}=\frac{\sqrt{p^{n-1}}}{n-1}\times\frac{a_{E/\F_p,n}}{2\sqrt{p^n}}+\frac{1}{2}\Big(\sqrt{p}+\frac{n-5}{(n-1)\sqrt{p}}\Big)
$$
Therefore, by Theorem\,\ref{thm6}, we conclude that $\Delta_{E/\F_p,n}$ is asymptotically given by
$$
\frac{\sqrt{p^{n-1}}}{n-1}\times\frac{(5-n)+(n-1)a_{E/\F_p,1}-(n-1)p}{2\sqrt{p^{n}}}+\frac{1}{2}\Big(\sqrt{p}+\frac{n-5}{(n-1)\sqrt{p}}\Big)$$
which is nothing but $\dis{\frac{a_{E/\F_p,1}}{2\sqrt{p}}}$, by a routine but direct simplification.
That is to say, we have proved the following:
\begin{thm}[Structure of  Big $\Delta$ in the rank $n$]\label{thm9}
With the same notation as above, we have for a fixed $n$, when $p$ becomes sufficiently large, we have asymptotically
$$\Delta_{E/\F_p,n}=\frac{a_{E/\F_p,1}}{2\sqrt{p}}.$$
\end{thm}
With this, it is now crystal clear that the rank $n$ Sato-Tate conjecture is a direct consequence of the classical Sato-Tate, namely, Theorem\,\ref{thm9}.\hskip 2.0cm $\square$~\\

\noindent
{\bf Acknowledgement}: The junior author was partially supported by the WISE program (MEXT) at Kyushu University. Primitive data of elliptic curve used here are from LMFDB,
while all programs are built under SageMath.

\vskip 16.0cm
Zhan SHI,\hskip 6.350cm Lin WENG,\\
shi.zhan.655@s.kyushu-u.ac.jp\hskip 3.30cm weng@math.kyushu-u.ac.jp\\
Graduate Program of Mathematics for Innovation\ \ Faculty of Mathematics\\ 
Kyushu University\hskip 5.050cm Kyushu University\\
Fukuoka, Japan\hskip 5.450cm Fukuoka, Japan
\end{document}